\documentclass[11pt]{amsart}
\usepackage{todonotes}

\usepackage{fullpage,cite,enumitem}
\usepackage{amssymb,amsmath,amsthm,mathtools,extpfeil}
\usepackage[all,cmtip]{xy}
\usepackage{rotating,graphicx}
\usepackage{setspace}

\usepackage[normalem]{ulem}

\usepackage{adjustbox}

\usepackage{tikz}
\usetikzlibrary{cd}
\usetikzlibrary{decorations.markings}
\usetikzlibrary{positioning}

\usepackage[bookmarks]{hyperref}
\hypersetup{
  colorlinks=true,
  linkcolor=blue,
  citecolor=blue,
  urlcolor=blue
}

\newtheorem{theorem}{Theorem}[section]
\newtheorem*{theorem*}{Theorem}
\newtheorem{theoremA}{Theorem}

\newtheorem{lemma}[theorem]{Lemma}
\newtheorem{conjecture}[theorem]{Conjecture}
\newtheorem{proposition}[theorem]{Proposition}

\theoremstyle{definition}

\numberwithin{equation}{section}

\newcommand{\Q}{\mathbb Q}

\newcommand{\lk}{\operatorname{lk}}

\begin{document}

\title{
The minimum size of the 3-dimensional
torus
}



\author{Seunghun Lee}
\address{Department of Mathematics, Keimyung University, Daegu, Republic of Korea}
\curraddr{}
\email{seunghun.math@gmail.com}

\author{Geunho Lim}
\address{Kangnam University, 40 Gangnam-ro, Giheung-gu, Yongin-si, Gyeonggi-do, 16979, Republic of Korea}
\curraddr{}
\email{limg@kangnam.ac.kr}

\author{Eran Nevo}
\address{The Hebrew University of Jerusalem,
Jerusalem, Israel and Institute of Mathematics, Universidad de Valladolid, Valladolid, Spain}
\email{nevo@math.huji.ac.il}

\dedicatory{Dedicated to the memory of Frank Hagen Lutz}

\begin{abstract}
We prove a conjecture of K\"uhnel, Lassmann, and Lutz, first asked in 1984, that the minimal number of vertices in a triangulation of the $3$-dimensional torus is $15$.
\end{abstract}

\setcounter{tocdepth}{1}

\maketitle
\setstretch{1.1}

\section{Introduction and main result}
Determining the minimal number of vertices in a simplicial triangulation of a given manifold is a classical problem in combinatorial topology~\cite{Brehm-Kuhnel:1987,Lutz:2005,datta:2007}.
Lutz, Sulanke, and Swartz~\cite{Lutz-Sulanke-Swartz:2009} 
determined this minimum for several $3$-manifolds, and Bagchi and Datta~\cite{BagchiDatta2008} for some sphere bundles.
More recent work on real projective spaces has produced subexponential constructions and superpolynomial vertex lower bounds~\cite{Adiprasito-Avvakumov-Karasev:2022,Frick-et-al:2026}, exact minimum up to dimension 4~\cite{Walkup1970, Lutz:2005, Balagopalan2017} and a conjectured minimum in dimension 5~\cite{GuyerSteinerbergerYang2026}.
See Lutz's \emph{Manifold Page}~\cite{Lutz:ManifoldPage} for tables of triangulations and further examples.

Nevertheless, even in dimensions three, many exact minima remain unknown, including the 3-torus.
For the two-dimensional torus, it is well-known that the minimal number is seven, obtained by the Cs\'asz\'ar torus~\cite{Csaszar:1949}.
The triangulation is the unique minimal one up to simplicial isomorphism~\cite[Section~3]{Lutz:2005}; see Figure~\ref{Figure-T7}.
In 1984, K\"uhnel and Lassmann constructed a triangulation of $T^3=S^1\times S^1\times S^1$ with $15$ vertices and $f$-vector $(15,105,180,90)$, and asked whether $15$ is the minimal number of vertices \cite{Kuehnel-Lassmann:1984}.
Lutz conjectured the minimality of this construction, and its uniqueness for its f-vector, in 2005 \cite[Conjectures~21 and~22]{Lutz:2005}.
We prove its minimality, leaving the uniqueness question open:
\begin{theoremA}\label{thm:minimum}
The minimal number of vertices in a triangulation of the $3$-dimensional torus is $15$.
\end{theoremA}

The construction of K\"uhnel and Lassmann provides the upper bound.
Our contribution is to rule out a triangulation with $14$ vertices.
This suffices since the stellar subdivision of a tetrahedron introduces one new vertex without changing the underlying manifold, so a triangulation with fewer than $14$ vertices would yield one with exactly $14$ vertices.

If $K$ is a triangulation of $T^3$ with $n$ vertices and $m$ edges, the lower bound inequality of Novik and Swartz~\cite[Theorem~5.2]{Novik-Swartz:2009} gives
\begin{equation}\label{eq:Novik-Swartz}
g_2(K)=m-4n+10\ge 10\beta_1(T^3)=30.
\end{equation}
Together with $m\le\binom n2$, this only implies $n\ge13$.
For $n=14$, it leaves a range of possible edge numbers.
Our proof uses the topology of induced subcomplexes to eliminate the $n=14$ possibility.
It cannot be done by a brute-force search of complexes on 14 vertices, since it is beyond reach of current computational power.
Here is an outline of our proof.

\textbf{Proof sketch.}
Partition the vertices of a hypothetical $14$-vertex triangulation $K$ into two sets of seven, and write $A$ and $B$ for the corresponding induced subcomplexes.
One can show that, if the inclusion $A\hookrightarrow T^3$ induces a nonzero map on $H_2(-;\Q)$, then $A$ must be the seven-vertex torus.
Duality for complementary induced subcomplexes then forces $B$ to be another seven-vertex torus.
We show that two such tori cannot have complementary vertex sets in a triangulated closed $3$-manifold.
Consequently, for every seven-vertex set $U$, the map
\[
H_2(K[U];\Q)\longrightarrow H_2(T^3;\Q)
\]
vanishes.
This vanishing map implies, by using standard algebraic topological techniques, that $\beta_1(K[U])\ge3$ and  $\beta_2(K[U])=c(K[[15]\setminus U])-1$, where $c(X)$ denotes the number of connected components of $X$.
Thus we obtain an inequality for the Euler characteristics of the induced subcomplexes.
Summing the resulting Euler characteristic inequalities over all $\binom{14}{7}$ choices of $U$ gives a lower bound on the ``total disconnectedness" $D$ of seven-vertex induced subgraphs.
More precisely, let $q=\binom{15}{2}-|K^{(1)}|$ be the number of missing edges of $K$, and let
\[
D:=\sum_{\substack{U\subseteq V(K)\\ |U|=7}}
   \bigl(c(K
   [U])-1\bigr).
\]
Summing up each side of the inequality for $\chi(K[U])$ over all $U$ and counting the contribution of each simplex in $K$ to this sum gives $D\ge198+126q$.
On the other hand, counting components of size one, two, and three in the seven-vertex induced subgraphs gives respective upper bounds of $20q$, $24q$, and $20q$.
Hence $D\le64q$, contradicting $q\ge0$.

\subsection*{Organization of the paper}
Section~2 establishes the seven-vertex input, the duality statement for complementary induced subcomplexes, the obstruction to two complementary seven-vertex tori, and the graph-counting estimate.
Section~3 combines these ingredients to exclude a $14$-vertex triangulation and concludes the proof of Theorem~\ref{thm:minimum}. Section~4 concludes with an approach to the minimality conjecture for higher dimensional tori.

\subsection*{Use of AI}
The authors used ChatGPT-6 Sol to assist with typesetting, language editing of portions of the text, proofreading, and locating references.
The authors take full responsibility for all content of this paper.

\section{Small induced subcomplexes and duality}

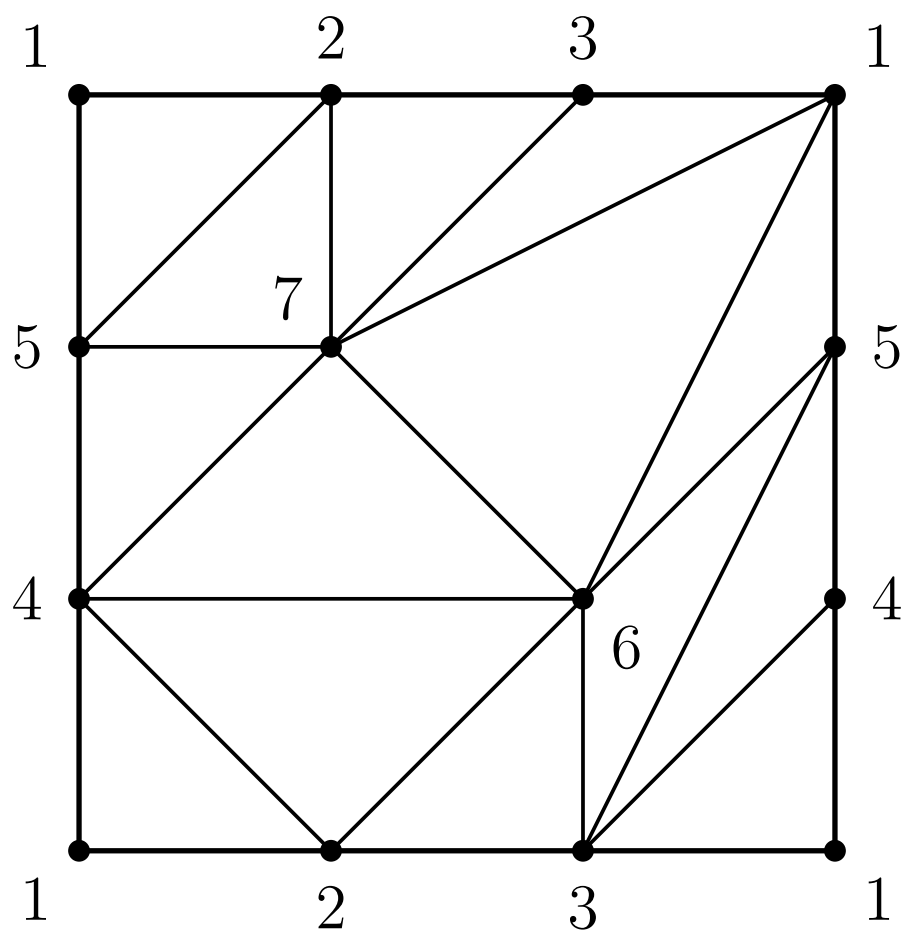
\begin{figure}
    \centering
\begin{tikzpicture}[
  x=5cm,y=5cm,
  line cap=round,line join=round,
  every node/.style={font=\large},
  edge/.style={line width=0.7pt},
  boundary/.style={line width=1pt},
  vertex/.style={circle,fill=black,inner sep=1.45pt,outer sep=0pt}
]
  \coordinate (bl) at (0,0);
  \coordinate (b2) at ({1/3},0);
  \coordinate (b3) at ({2/3},0);
  \coordinate (br) at (1,0);
  \coordinate (l4) at (0,{1/3});
  \coordinate (l5) at (0,{2/3});
  \coordinate (tl) at (0,1);
  \coordinate (t2) at ({1/3},1);
  \coordinate (t3) at ({2/3},1);
  \coordinate (tr) at (1,1);
  \coordinate (r4) at (1,{1/3});
  \coordinate (r5) at (1,{2/3});

  \coordinate (v6) at ({2/3},{1/3});
  \coordinate (v7) at ({1/3},{2/3});

  \foreach \u/\v in {
    b2/l4,b2/v6,b3/r4,b3/r5,b3/v6,
    l4/v6,l4/v7,l5/t2,l5/v7,r5/v6,
    t2/v7,t3/v7,tr/v6,tr/v7,v6/v7%
  }{\draw[edge] (\u)--(\v);}

  \draw[boundary] (bl)--(b2)--(b3)--(br)--(r4)--(r5)--(tr)
                  --(t3)--(t2)--(tl)--(l5)--(l4)--cycle;

  \foreach \p in {bl,b2,b3,br,l4,l5,tl,t2,t3,tr,r4,r5,v6,v7}
    \node[vertex] at (\p) {};

  \node[below left=2pt]  at (bl) {$1$};
  \node[below right=2pt] at (br) {$1$};
  \node[above left=2pt]  at (tl) {$1$};
  \node[above right=2pt] at (tr) {$1$};

  \node[below=3pt] at (b2) {$2$};
  \node[below=3pt] at (b3) {$3$};
  \node[above=3pt] at (t2) {$2$};
  \node[above=3pt] at (t3) {$3$};
  \node[left=3pt]  at (l4) {$4$};
  \node[left=3pt]  at (l5) {$5$};
  \node[right=3pt] at (r4) {$4$};
  \node[right=3pt] at (r5) {$5$};

  \node[below right=2pt] at (v6) {$6$};
  \node[above left=2pt]  at (v7) {$7$};

\end{tikzpicture}
\caption{The Cs\'asz\'ar torus $T_7$.}
    \label{Figure-T7}
\end{figure}

All homology and cohomology groups in this paper have coefficients in $\Q$.
Write $\beta_i(X)=\dim_\Q H_i(X;\Q)$ and $c(X)=\beta_0(X)$.
For a simplicial complex $K$ and a vertex subset $U$, write $K[U]$ for the induced subcomplex on $U$.
Let $T_7$ denote the Cs\'asz\'ar triangulation of the torus; see Figure~\ref{Figure-T7}.
 
We will make use only of the following two facts on $T_7$: 
its $f$-vector is $(7,21,14)$, 
and every vertex link of it is a six-cycle.

\subsection{A seven-vertex rigidity result}

We show that a nonzero image in the second homology of an aspherical space forces a complex on at most seven vertices to be the seven-vertex torus.

\begin{lemma}\label{lem:seven}
Let $X$ be a simplicial complex on at most seven vertices embedded in $T^3$. If the embedding $\iota: X\to T^3$ induces a non-zero homomorphism $\iota_*:H_2(X)\to H_2(T^3)$, then $X$ is isomorphic to $T_7$.
\end{lemma}

\begin{proof}
 Let $\pi=\pi_1(X)$ and $B\pi$ denote the classifying space of $\pi$.
We rely on some technical but well-known facts from algebraic topology for connected $X$.
The first is that the map $\iota_*$ factors through $H_2(B\pi)$, since $T^3$ is aspherical~\cite[Example~1B.5]{Hatcher:2002}.
(For more details, see~\cite[Theorem~7.26]{Davis-Kirk:2001}.)
In other words, there are maps $c_X:X\to B\pi_1(X)$ and $\overline{f}:B\pi_1(X)\to M$ such that the following diagram commutes up to homotopy:
\[
\begin{tikzcd}[column sep=large, row sep=large]
X \arrow[r, "c_X"] \arrow[dr, "\iota"'] & B\pi \arrow[d, "\bar f"] \\
& T^3
\end{tikzcd}
\]
The second is that when $\pi$ is a free product of a cyclic group and a free group, we have $H_2(B\pi)=0$ (See~\cite[Theorem~1B.11,  Corollary~2.25, Example~1B.1, and p.~148]{Hatcher:2002}).

Now, if the fundamental group of a connected simplicial complex $X$ on at most 7 vertices is not of this product type, then $X$ must be $T_7$ by \cite[Theorem 3.10]{Govc-Marzantowicz-Michalak-Pavesic:2025}.
Combining all these facts and considering a connected component of $X$ if necessary readily yields the lemma.
\end{proof}

\subsection{Complementary induced subcomplexes}

We use Poincar\'e--Lefschetz duality to obtain constraints on the Betti numbers of complementary induced subcomplexes when their inclusions into the ambient manifold induce zero maps on second homology.

\begin{lemma}\label{lem:duality}
Let $K$ triangulate a connected closed orientable $3$-manifold $M$, let $\varnothing\ne U\subsetneq V(K)$, and put
\[
A=K[U],\qquad B=K[V(K)\setminus U].
\]
Then there is an exact sequence
\[
0=H_3(B)\longrightarrow H_3(M)\longrightarrow H^0(A)\longrightarrow H_2(B)
\longrightarrow H_2(M)\longrightarrow H^1(A).
\]
If both maps $H_2(A)\to H_2(M)$ and $H_2(B)\to H_2(M)$ vanish, then
\[
\beta_1(A)\ge \beta_2(M)\qquad \text{and} \qquad
\beta_2(A)=c(B)-1.
\]
\end{lemma}

\begin{proof}
Since $A$ is an induced subcomplex of $K$, the space $|B|$ is a strong deformation retract of $M\setminus |A|$ by~\cite[Lemma~70.1]{Munkres:1984}.
To see this, write a point outside $|A|$ in barycentric coordinates and continuously scale the coordinates on $U$ to zero, renormalizing the other coordinates.
The sum of the coordinates outside $U$ is positive because $A$ is induced.
The homotopy stays within the simplex containing the point and fixes $|B|$.

Let $N$ be a closed regular neighborhood of $|A|$.
Excision and Poincar\'e--Lefschetz duality give
\[
H_i(M,M\setminus |A|)
\cong H_i(N,\partial N)
\cong H^{3-i}(N)
\cong H^{3-i}(A).
\]
Plugging into the long exact sequence of the pair $(M, M\setminus |A|)$  the isomorphisms $H_i(M\setminus |A|)\cong H_i(B)$ and $H_i(M, M\setminus |A|)\cong H^{3-i}(A)$
gives the displayed sequence in the Lemma.
Since $B$ is a proper subcomplex of a triangulation of $M$, we have $H_3(B)=0$.

If $H_2(B)\to H_2(M)$ is zero, exactness gives an injection $H_2(M)\hookrightarrow H^1(A)$.
By the universal coefficient theorem, we obtain $\beta_1(A)\ge \beta_2(M)$.
Interchanging $A$ and $B$, the vanishing of $H_2(A)\to H_2(M)$ gives $H_2(A)\cong H^0(B)/\Q$, where $\Q$ is the diagonal subspace.
Applying the universal coefficient theorem again, we obtain $
\beta_2(A)=\beta_0(B)-1=c(B)-1$.
\end{proof}

\subsection{Excluding two complementary seven-vertex tori}

We show that a nonzero ambient second-homology class on seven vertices would force two complementary copies of $T_7$, which vertex-link counting rules out.

\begin{lemma}\label{lem:patch}
If $C$ is a connected contractible subcomplex of $T_7$ on $r$ vertices, where $2\le r\le6$, then the number of triangles in $C$ is bounded from above by
\[
f_2(C)\le r-2.
\]
\end{lemma}

\begin{proof}
For every vertex $x$ of $C$, the link $\lk_C(x)$ is a proper subcomplex of the six-cycle $\lk_{T_7}(x)$.
It is therefore a forest.
It is nonempty because $C$ is connected and $r\ge2$.
Thus
\begin{equation}\label{eq:link-bound}
f_1(\lk_C(x))\le f_0(\lk_C(x))-1.
\end{equation}
Each $2$-simplex of $C$ contributes one edge to the link of each of its three vertices, so
\[
\sum_{x\in V(C)} f_1(\lk_C(x))=3f_2(C).
\]
Each edge of $C$ contributes one vertex to the link of each endpoint, so
\[
\sum_{x\in V(C)} f_0(\lk_C(x))=2f_1(C).
\]
Summing up each side of the inequality~(\ref{eq:link-bound}) over all $r$ vertices yields
\[
3f_2(C)\le 2f_1(C)-r.
\]
Since $C$ is contractible, $r-f_1(C)+f_2(C)=1$.
Substitution gives $f_2(C)\le r-2$.
\end{proof}

\begin{lemma}\label{lem:two-tori}
A triangulated closed $3$-manifold $K$ on fourteen vertices cannot have its fourteen vertices partitioned into sets $U,W$ such that both $K[U]$ and $K[W]$ are isomorphic to $T_7$.
\end{lemma}

\begin{proof}
Suppose for contradiction that such a partition exists, and write $A=K[U]$ and $B=K[W]$.
Let $a$ be a vertex in $U$.
The induced subcomplex $\lk_K(a)[U\setminus\{a\}]$ is the six-cycle $\lk_A(a)$.
By the Jordan curve theorem, its complement in $|\lk_K(a)|\cong S^2$ has two components, each is an open disk.
The induced-complement deformation retraction argument~\cite[Lemma~70.1]{Munkres:1984} used in Lemma~\ref{lem:duality}, applied to the triangulation $\lk_K(a)$ of the $2$-sphere and its proper nonempty vertex subset $U\setminus \{a\}$, shows that $\lk_K(a)[W]$ has exactly two connected components, both are contractible.
Each component has at most six vertices because $\lk_K(a)[W]$ has at most seven vertices in total.

The local picture at $a$ is illustrated in Figure~\ref{fig:two-tori-link}.
\begin{figure}[htbp]
\centering
\begin{tikzpicture}[
  x=1cm,y=1cm,font=\small,
  every node/.style={inner sep=2pt},
  vertex/.style={circle,fill=black,inner sep=1.5pt},
  wvertex/.style={circle,fill=teal!70!black,inner sep=1.7pt}
]

\filldraw[fill=orange!10,draw=black!65,thick]
  (-4.4,-1.65)--(2.6,-1.65)
  --(4.4,1.45)--(-2.6,1.45)--cycle;

\fill[blue!4] (0,0) circle[radius=2];

\begin{scope}
  \clip (0,0) circle[radius=2];
  \fill[orange!10,opacity=0.75]
    (-2,0)
    arc[start angle=180,end angle=360,
        x radius=2,y radius=0.6]
    --(2,-1.65)--(-2,-1.65)--cycle;
\end{scope}

\fill[orange!16]
  (0,0) ellipse[x radius=2,y radius=0.6];

\foreach \i/\angle in {1/0,2/60,3/120,4/180,5/240,6/300}
  \coordinate (u\i) at
    ({2*cos(\angle)},{0.6*sin(\angle)});

\coordinate (a) at (0,0);

\foreach \i in {2,3}
  \draw[black!30,densely dashed] (a)--(u\i);
\foreach \i in {1,4,5,6}
  \draw[black!30] (a)--(u\i);

\pgfmathsetmacro{\frontangle}{asin(1.65/2)}

\draw[thick] (2,0)
  arc[start angle=0,end angle=180,radius=2];

\draw[thick,dashed,black!65] (-2,0)
  arc[start angle=180,
      end angle={180+\frontangle},radius=2];

\draw[thick]
  ({2*cos(180+\frontangle)},-1.65)
  arc[start angle={180+\frontangle},
      end angle={360-\frontangle},radius=2];

\draw[thick,dashed,black!65]
  ({2*cos(360-\frontangle)},-1.65)
  arc[start angle={360-\frontangle},
      end angle=360,radius=2];

\draw[black!65,thick]
  (-4.4,-1.65)--(2.6,-1.65);

\draw[orange!80!black,thick,densely dashed] (2,0)
  arc[start angle=0,end angle=180,
      x radius=2,y radius=0.6];

\draw[orange!80!black,thick] (-2,0)
  arc[start angle=180,end angle=360,
      x radius=2,y radius=0.6];

\foreach \i in {1,...,6}
  \node[vertex] at (u\i) {};

\node[vertex,label=below right:$a$] at (a) {};

\coordinate (w1) at (-0.85,0.95);
\coordinate (w2) at (0.15,0.78);
\coordinate (w3) at (0.85,1.10);
\coordinate (w4) at (-0.15,1.62);

\filldraw[fill=teal!15,draw=teal!70!black,thick]
  (w1)--(w2)--(w3)--(w4)--cycle;
\draw[teal!70!black,thick] (w1)--(w3);

\foreach \i in {1,...,4}
  \node[wvertex] at (w\i) {};

\coordinate (w5) at (-0.60,-1.25);
\coordinate (w6) at (0.60,-1.45);
\draw[teal!70!black,thick] (w5)--(w6);
\node[wvertex] at (w5) {};
\node[wvertex] at (w6) {};

\node at (-0.25,2.40)
  {$|\lk_K(a)|\cong S^2$};

\node[anchor=east] at (-2.35,-0.20)
  {$\lk_A(a)$};
\draw[->,orange!80!black,thick]
  (-2.30,-0.20)--(-1.72,-0.30);

\node[anchor=west,align=left] at (2.65,2.05)
  {two components of\\$\lk_K(a)[W]$};
\draw[teal!70!black,thick]
  (2.60,2.05)--(2.25,1.75);
\draw[->,teal!70!black,thick]
  (2.25,1.75)--(0.96,1.15);
\draw[->,teal!70!black,thick] (2.25,1.75)
  .. controls (2.95,-1.60) and (1.65,-1.68)
  .. (0.74,-1.47);

\node[anchor=west] at (2.50,-2.10)
  {local patch of $|A|$};
\draw[->,orange!80!black,thick]
  (2.65,-1.85)--(2.70,-1.35);

\end{tikzpicture}
\caption{The local picture at $a\in U$.
The link sphere meets the local surface patch of $|A|$ along the induced six-cycle $\lk_A(a)$.
The surface patch is drawn translucently to show both contractible components of $\lk_K(a)[W]$.}
\label{fig:two-tori-link}
\end{figure}
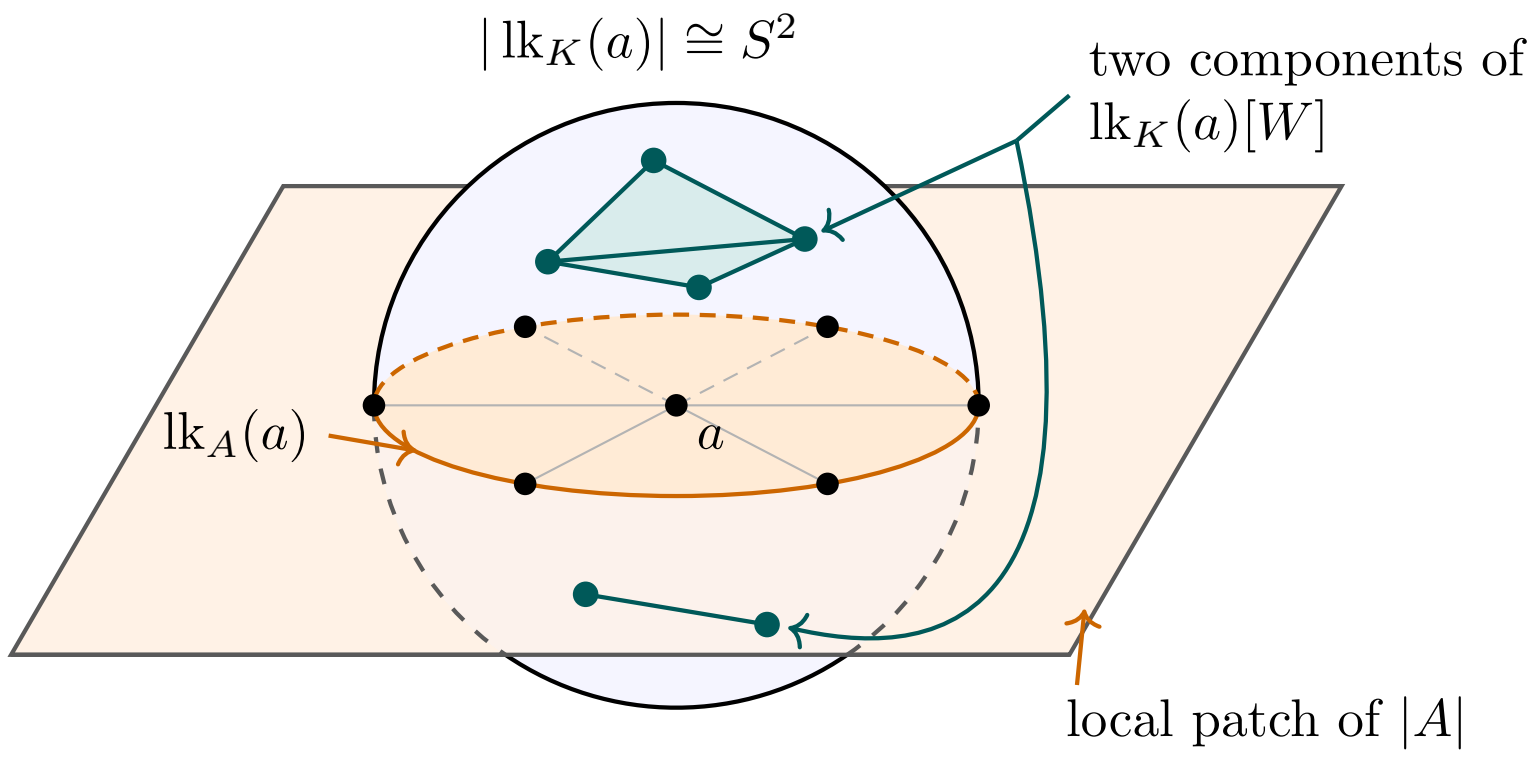

Neither component of $\lk_K(a)[W]$ is a single vertex.
Indeed, suppose one component is the singleton $b\in W$.
Then the triangulation of the closed disk in $\lk_K(a)$ containing $b$ with boundary $\lk_A(a)$ must be the cone $\{b\}*\lk_A(a)$, as $\lk_A(a)$ is an induced subcomplex of $\lk_K(a)$.
It follows that $abxy$ is a tetrahedron of $K$ for every edge $xy$ of $\lk_A(a)$.
Consequently, $\lk_K(b)[U]$ contains the entire closed star of $a$ in $A$.
This star is connected and contains all seven vertices of $U$.
But the preceding argument with $U$ and $W$ interchanged says that $\lk_K(b)[U]$ has two nonempty components, a contradiction.

Let $r_a=f_0(\lk_K(a)[W])$.
The last paragraph implies that each of the two components of $\lk_K(a)[W]$ satisfies the conditions of Lemma~\ref{lem:patch}.
Applying this lemma and summing over the two components gives
\[
f_2(\lk_K(a)[W])\le r_a-4. 
\]
Let $e$ denote the number of edges of $K$ with one vertex in $U$ and the other vertex in $W$.
Then
\[
\sum_{a\in U}r_a=e\le 7\cdot 7=49.
\]
Every triangle of $B=K[W]$ lies in exactly two tetrahedra of $K$, and their fourth vertices belong to $U$ because $B\cong T_7$ is two-dimensional.
Since $f_2(B)=14$, we obtain
\[
28=14\cdot 2=\sum_{a\in U}f_2(\lk_K(a)[W])
\le\sum_{a\in U}(r_a-4)
=e-28
\le 49-28=  21.
\]
This contradiction proves the lemma.
\end{proof}

\begin{proposition}\label{prop:vanishing}
If $K$ triangulates $T^3$ on fourteen vertices, then for every seven-element subset $U$ the induced map
\[
H_2(K[U])\longrightarrow H_2(T^3)
\]
is zero.
For $A=K[U]$ and $B=K[V(K)\setminus U]$, one consequently has
\[
\beta_1(A)\ge3,\qquad
\beta_2(A)=c(B)-1.
\]
\end{proposition}

\begin{proof}
Suppose the above map on $H_2$ is nonzero for some $U$.
By Lemma~\ref{lem:seven}, $A=K[U]$ is $T_7$.
In particular, $\dim H^1(A)=2$.
Since $\dim H_2(T^3)=3$, the map $H_2(T^3)\to H^1(A)$ in Lemma~\ref{lem:duality} has a nonzero kernel.
Exactness forces $H_2(B)\to H_2(T^3)$ to have a nonzero image.
Lemma~\ref{lem:seven} then implies that $B$ is also $T_7$, contradicting Lemma~\ref{lem:two-tori}.
The final assertions of this lemma follow from the final assertions of Lemma~\ref{lem:duality}.
\end{proof}

\subsection{Counting disconnected induced subgraphs}

We bound the total number of components beyond the first in seven-vertex induced subgraphs in terms of the number of missing edges.
As we shall detail in the next section, the 1-skeleton of a hypothetical 14-vertex triangulation of $T^3$ must satisfy the conditions in the lemma below.
\begin{lemma}\label{lem:graph}
Let $G$ be a simple graph on fourteen vertices with minimum degree at least four, and suppose its complement $G^c$ has $q\le15$ edges.
Define
\[
D=\sum_{\substack{U\subseteq V(G)\\|U|=7}}
\bigl(c(G[U])-1\bigr).
\]
Then $D\le64q$.
\end{lemma}

\begin{proof}
In a disconnected graph on seven vertices, every component other than a chosen largest component has at most three vertices.
For $k=1,2,3$, let $N_k$ count all occurrences of a component of size $k$ in the graphs $G[U]$, over all $U$ with $|U|=7$.
Thus $D\le N_1+N_2+N_3$.

First we bound $N_1$.
For $v\in V(G^c)$, write $d_v=\deg_{G^c}(v)$.
The minimum-degree assumption on $G$ gives $d_v\le9$, and $\sum_v d_v=2q$.
Then, as $N_1$ is the number of pairs $(v,U)$ where $|U|=7$ and $v$ is an isolated vertex in $G[U]$, we get
\[
N_1=\sum_v\binom{d_v}{6}\le10\sum_vd_v=20q,
\]
where the inequality follows from $\binom d6\le10d$ for each $0\le d\le9$.

Next we bound $N_2$.
For an unordered pair $\{x,y\}$, let $t_{xy}=|N_{G^c}(x)\cap N_{G^c}(y)|$.
The $2t_{xy}$ edges joining the pair to its common neighbors in $G^c$ are distinct, so $t_{xy}\le\lfloor q/2\rfloor\le  \lfloor 15/2\rfloor=7$.
Note that $\sum_{\{x,y\}}t_{xy}=\sum_v\binom{d_v}{2}$, since both sides count pairs $(v,\{x,y\})$ in which $x$ and $y$ are distinct neighbors of $v$ in $G^c$.
Dropping the requirement that the pair be connected in $G$ only increases the count, and therefore
\[
\begin{aligned}
N_2
&\le\sum_{\{x,y\}}\binom{t_{xy}}5
\le3\sum_{\{x,y\}}t_{xy}\\
&=3\sum_v\binom{d_v}2
\le12\sum_vd_v
=24q.
\end{aligned}
\]
Here $\binom t5\le3t$ for $0\le t\le7$, and $\binom d2\le4d$ for $d\le9$.

Similarly, we bound $N_3$.
For an unordered triple $S$, put $t_S=|\bigcap_{x\in S}N_{G^c}(x)|$.
The $3t_S$ edges between $S$ and its common neighbors are distinct, so $t_S\le\lfloor q/3\rfloor\le5$.
Also, note that $\sum_{|S|=3}t_S=\sum_v\binom{d_v}{3}$, since both sides count pairs $(v,S)$ in which $S$ is a three-element subset of the neighbors of $v$ in $G^c$.
It follows that
\[
\begin{aligned}
N_3
&\le\sum_{|S|=3}\binom{t_S}4
\le\sum_{|S|=3}t_S\\
&=\sum_v\binom{d_v}3
\le10\sum_vd_v
=20q.
\end{aligned}
\]
For the two inequalities we used $\binom t4\le t$ for $0\le t\le5$, and $\binom d3\le10d$ for $d\le9$, resp.
Adding these three estimates gives $D\le64q$.
\end{proof}

\section{Proof of Theorem~\ref{thm:minimum}}

In this section, we prove Theorem~\ref{thm:minimum}
by assuming that $T^3$ has a fourteen-vertex triangulation $K$ and deriving a contradiction.

\begin{theorem}\label{thm:main}
There is no simplicial triangulation of $T^3$ on fourteen vertices.
\end{theorem}

\begin{proof}
Suppose that $K$ is such a triangulation, and set $m=f_1(K)$.
The 
Klee-Dehn–Sommerville relations~\cite{Klee-DS},
first proved by Poincar\'e 
for closed $3$-manifolds~\cite{Poincare-AnalysisSitus}, give
\[
f(K)=(14,m,2m-28,m-14).
\]
The Novik-Swartz lower bound inequality (\ref{eq:Novik-Swartz}) gives
\[
m-46=g_2(K)\ge10\beta_1(T^3)=30.
\]
Hence
$76\le m\le91=\binom{14}{2}$, 
thus the number of missing edges satisfies
$q:=91-m\in\{0,\ldots,15\}$.
Let $G=K^{(1)}$ be the 1-skeleton and $D$ as defined in Lemma~\ref{lem:graph}.
Every vertex link is a simplicial two-sphere and thus has at least four vertices, so the minimum degree of $G$ is at least four.
Consequently, by Lemma~\ref{lem:graph}, 
\begin{equation}\label{eq:D}
D\le64q.
\end{equation}
For every seven-element vertex subset $U$, set
\[
A_U=K[U],\qquad B_U=K[V(K)\setminus U].
\]
Proposition~\ref{prop:vanishing} gives $\beta_1(A_U)\ge3$ and $\beta_2(A_U)=c(B_U)-1$.
Moreover, $\beta_3(A_U)=0$ because $A_U$ is a proper subcomplex of a connected closed triangulated $3$-manifold.
Indeed, a nonzero top-dimensional cycle in such a manifold has nonzero coefficients on every tetrahedron.
Therefore,
\begin{equation}\label{eq:xi-ineq}
\chi(A_U)
=c(A_U)-\beta_1(A_U)+\beta_2(A_U)
\le c(A_U)+c(B_U)-4.
\end{equation}
There are $N=\binom{14}{7}
$ seven-vertex subsets, and the involution $U\mapsto V(K)\setminus U$ 
permutes them.
Summing over $U$ each side of the inequality (\ref{eq:xi-ineq}) gives \begin{equation}\label{eq:2D-2N}
\sum_{|U|=7}\chi(A_U)\le 2D-2N.
\end{equation}
Each $i$-simplex of $K$ belongs to exactly $\binom{13-i}{6-i}$ of these seven-vertex induced subcomplexes.
Consequently,
\[
\begin{aligned}
\sum_{|U|=7}\chi(A_U)
&=14\binom{13}{6}
-m\binom{12}{5}
+(2m-28)\binom{11}{4}
-(m-14)\binom{10}{3}\\
&=16464-252m.
\end{aligned}
\]
Substituting into (\ref{eq:2D-2N}) yields a lower bound on $D$:
\[
D\ge11664-126m= 11664-126 (91-q)=  198+126q.
\]
Combining it with the upper bound on $D$, namely (\ref{eq:D}), we obtain
$198+126q\le64q$.
This is a contradiction to the assumption that $q \ge 0$.
\end{proof}

Theorem~\ref{thm:main}, together with the stellar subdivision observation in the introduction and the $15$-vertex construction of K\"uhnel and Lassmann \cite{Kuehnel-Lassmann:1984}, proves Theorem~\ref{thm:minimum}.

\section{Approach to minimal tori in higher dimensions}
K\"uhnel and Lassmann~\cite{Kuhnel-Lassmann-1988} constructed for every $d\ge 2$ a triangulation $K(T^d)$ of the $d$-dimensional torus $T^d$ with $2^{d+1}-1$ vertices, which has a complete 1-skeleton. Lutz~\cite[Conjecture 21]{Lutz:2005} conjectured that $K(T^d)$ is vertex minimal.
As an approach to prove this conjecture, potentially by following the strategy in our proof of the $d=3$ case, we propose the following.
\begin{conjecture}\label{conj:minimal-T^d}
Suppose that $K$ is a simplicial complex whose cohomology ring is isomorphic to the cohomology ring of $T^d$.
Then $K$ has at least $2^{d+1}-1$ vertices.
\end{conjecture}
Note that uniqueness of the vertex minimal triangulation fails in dimension four, as shown by Grigis~\cite{Grigis1998}.

\bibliographystyle{amsalpha-order}
\bibliography{reference}{}
\def\MR#1{}
\end{document}